# Об эффективных рандомизированных алгоритмах поиска вектора PageRank[1]

*А.В. Гасников, Д.Ю. Дмитриев*

gasnikov@yandex.ru, dmitden@gmail.com

Лаборатория структурных методов анализа данных в предсказательном моделировании,

Кафедра математических основ управления,

Факультет управления и прикладной математики МФТИ

Институт проблем передачи информации РАН

В работе рассматриваются два рандомизированных способа поиска вектора PageRank, т.е. решения системы $\vec{p}^T = \vec{p}^T P$, со стохастической матрицей $P$ размера $n \times n$ (решение ищется в классе распределений вероятностей), где $n \sim 10^7 - 10^9$, с точностью $\varepsilon$: $\varepsilon \gg n^{-1}$. Таким образом, исключается возможность "честного" умножения матрицы $P$ на столбец, если рассматривать не разреженные объекты. Первый способ базируется на идее Markov Chain Monte Carlo. Этот подход эффективен в случае "быстрого" выхода итерационного процесса $\vec{p}_{t+1}^T = \vec{p}_t^T P$ на стационар, и учитывает также другую специфику матрицы $P$ – равенство отличных от нуля вне диагональных элементов матрицы $P$ по строчкам (это используется при организации случайного блуждания по графу с матрицей $P$). На основе современных неравенств концентрации меры в работе приводятся новые оценки времени работы такого метода, учитывающие специфику матрицы $P$. В основе второго способа идея – свести поиск ранжирующего вектора к поиску равновесия в антагонистической матричной игре:

$$\min_{\vec{p} \in S_n(1)} \max_{\vec{u} \in S_n(1)} \left\langle \vec{u}, \left(P^T - I\right) \vec{p} \right\rangle,$$

где $S_n(1)$ – единичный симплекс в $\mathbb{R}^n$, а $I$ – единичная матрица. Возникшая задача решается с помощью небольшой модификации алгоритма Григориадиса–Хачияна (1995). Этот метод, также как метод Назина–Поляка (2009), является рандомизированным вариантом метода зеркального спуска А.С. Немировского. Отличие заключается в том, что у Григориадиса–Хачияна рандомизация осуществляется на этапе проектирования на симплекс, а не на этапе вычисления стохастического градиента. Для разреженных матриц $P$ предложенный нами метод показывает заметно лучшие результаты.

[1] Работа выполнена при финансовой поддержке РФФИ (код проекта 14-01-00722-а); Лаборатории структурных методов анализа данных в предсказательном моделировании ФУПМ МФТИ; гранта Президента РФ № МК-5285.2013.9. Исследования, связанные с получением оценок метода МСМС (см. п. 3) выполнены за счет гранта Российского научного фонда (проект № 14-50-00150).

# 1. Введение

Известно, что поисковая система Google была создана в качестве учебного проекта студентов Стэнфордского университета, см. [1]. В 1996 году авторы работали над поисковой системой BackRub, а в 1998 году на её основе создали новую поисковую систему Google [2], [3]. В [1] был предложен определенный способ ранжирования web-страниц. Этот способ, также как и довольно большой класс задач ранжирования, возникающих, например, при вычислении индексов цитирования ученых или журналов [4], сводится [5] к нахождению левого собственного вектора (нормированного на единицу: $\sum_{k=1}^{n} p_k = 1$), отвечающего собственному значению 1, некоторой стохастической (по строкам) матрицы $P = \| p_{ij} \|_{i,j=1}^{n,n}$: $\vec{p}^T = \vec{p}^T P$, $n \gg 1$

**Замечание 1.** К поиску такого вектора, который иногда называют вектором Фробениуса–Перрона, сводится (например, в модели де Гроота) задача поиска консенсуса. Подробнее об этом, и в целом о моделях консенсуса, написано в обзоре [6].

**Замечание 2.** Решение $\vec{p}^T = \vec{p}^T P$ всегда существует по теореме Брауера (непрерывный (ограниченный) оператор $P$ отображает выпуклый компакт (симплекс) в себя), и единственно в классе распределений вероятностей тогда и только тогда, когда имеется всего "один класс сообщающихся (существенных) состояний", при возможном наличии "несущественных состояний" [7]. Другими словами, если мы поставим в соответствие матрице $P$ такой ориентированный граф, что вершины $i$ и $j$ соединены ребром тогда и только тогда, когда $p_{ij} > 0$, то в таком графе любая вершина может принадлежать только одному из двух типов: несущественная – стартуя из этой вершины, двигаясь по ребрам с учетом их ориентации, всегда можно забрести в такую вершину, из которой обратно никогда не вернемся; существенная – стартуя из любой такой вершины, мы можем добраться в любую другую существенную вершину (в частности вернуться в исходную).

Приведем, следуя [5], обоснование такому способу. Пусть имеется ориентированный граф $G = \langle V, E \rangle$ сети Интернет (вершины – web-страницы, ребра – ссылки: запись $(i, j) \in E$ означает, что на $i$-й странице имеется ссылка на $j$-ю страницу), $N$ – число пользователей сети (это число не меняется со временем). Пусть $n_i(t)$ – число посетителей web-страницы $i$ в момент времени $t$. За один такт времени каждый посетитель этой web-страницы независимо ни от чего с вероятностью $p_{ij}$ переходит по ссылке на web-страницу $j$. Считаем стохастическую матрицу $P$ неразложимой и апериодической (см. [7]). Ниже приведен результат из [5] (обоснование будет приведено в п. 4), позволяющий по-другому интерпретировать вектор $\vec{p}$ (PageRank), согласно которому и происходит ранжирование web-страниц:

$$\exists\ \lambda_{0.99} > 0, T_G > 0:\ \forall\ t \ge T_G$$

$$P\left(\left\|\frac{\vec{n}(t)}{N} - \vec{p}\right\|_2 \le \frac{\lambda_{0.99}}{\sqrt{N}}\right) \ge 0.99,$$

где $\vec{p}^T = \vec{p}^T P$ (решение единственно в классе распределений вероятностей в виду неразложимости).

В ряде случаев считают, что

$$P = (1-\delta) I + \delta \tilde{P},$$

где $\delta \in (0,1]$, $I = \mathrm{diag}\{1,...,1\}$ – единичная матрица,

$$\tilde{p}_{ij} = \begin{cases} \left|\{k:\ (i,k) \in E\}\right|^{-1}, i \ne j \\ 0, \quad \text{иначе} \end{cases}.$$

Такая специфика матрицы $P$ нами будет использоваться в п. 4. Отметим также, что вместо $I$ часто берут стохастическую матрицу с одинаковыми элементами – матрицу телепортации (см. [2]). Это сразу дает оценку снизу $\delta$ на спектральную щель матрицы $P$.

Опишем вкратце структуру статьи. В п. 2 приводится обзор наиболее популярных известных ранее способов численного поиска вектора PageRank. Новизна заключается в том, что этот обзор делается одновременно с декларированием двух новых способов поиска вектора PageRank: на основе Markov chain Monte Carlo и на основе алгоритма Григориадиса–Хачияна. В п. 3, также носящем обзорный характер, кратко описываются основные необходимые в дальнейшем факты о методе Markov chain Monte Carlo (MCMC). Новых результатов в этом пункте нет. Тем не менее, интерес может представлять обзор литературы. В п. 4 метод MCMC применяется для поиска вектора PageRank. Новыми здесь являются оценки скорости сходимости такого метода, а также оценки общего числа затраченных арифметических операций. Подчеркнем, что в отличие от общего случая, мы ограничиваемся в данной статье изучением эффективности метода MCMC на специальном клас-се задач. За счет этого удается получить существенно лучшие оценки, чем можно было бы ожидать в общем случае. Новизна метода также заключается в том, как организуются случайные блуждания. Отметим, что предложена хорошо параллелизуемая версия метода MCMC. Описанный в предыдущем пункте метод MCMC будет хорошо работать, если спектральная щель матрицы $P$ достаточно велика. В п. 5 предложен новый способ поиска вектора PageRank, не требующий ограничений на спектральную щель. Этот способ сводит поиск ранжирующего вектора к поиску равновесия Нэша в антагонистической матричной игре. Для поиска равновесия мы используем метод Григориадиса–Хачияна, который позволяет учитывать разреженность матрицы $P$ и хорошо распараллеливается.

## 2. Обзор и обсуждение известных ранее и новых способов численного поиска вектора PageRank

В работах [8]–[16] предложены различные способы численного поиска вектора PageRank $\vec{p}_*$.

**Замечание 3.** Строго говоря, в [11] и [12] решались другие задачи, но несложно перенести алгоритмы этих работ на задачу поиска вектора PageRank: в случае [11] это делается тривиально [13], а вот в случае [12] потребовалось немного более точное исследование сходимости – поправка $\ln\left(\sigma^{-1}\right)$ появилась у нас из-за того, что мы избавились от математического ожидания в критерии качества (цели).

Приведем краткое резюме сложностных оценок алгоритмов работ [8]–[16] и алгоритмов, предложенных в этой статье. “Сложность” понимается как количество арифметических операций типа умножения двух чисел, которые достаточно осуществить, чтобы с вероятностью не меньше $1-\sigma$ достичь точности решения $\varepsilon$ по “Целевому” функционалу.

| **Метод** | **Условие** | **Сложность** | **Цель (min)** |
|---|---|---|---|
| Назина–Поляка [8] | нет | $O\left(\frac{n\ln\left(n/\sigma\right)}{\varepsilon^2}\right)$ | $\left\|P^T\vec{p}-\vec{p}\right\|_2^2$ |
| Ю.Е. Нестерова [9], [10], [13] | $\bar{S}$ | $O\left(\frac{s^2\ln n}{\varepsilon^2}\right)$ | $\left\|P^T\vec{p}-\vec{p}\right\|_2$ |
| Юдицкого–Лана–Немировского–Шапиро [11], [13] | нет | $O\left(\frac{n\ln(n/\sigma)}{\varepsilon^2}\right)$ | $\left\\|P^T\vec{p}-\vec{p}\right\\|_\infty$ |
| Григориадиса–Хачияна [12] | $\bar{S}$ | $O\left(\frac{s\ln n\ln(n/\sigma)}{\varepsilon^2}\right)$ | $\left\\|P^T\vec{p}-\vec{p}\right\\|_\infty$ |
| Нестерова– Немировского [14] | G, S | $\frac{sn}{\alpha}\ln\left(\frac{2}{\varepsilon}\right)$ | $\left\\|\vec{p}-\vec{p}_*\right\\|_1$ |
| Поляка–Трембы [15] | S | $\frac{2sn}{\varepsilon}$ | $\left\\|P^T\vec{p}-\vec{p}\right\\|_1$ |
| Д. Спилмана [16] | G, S | $O\left(\frac{s^2}{\alpha\varepsilon}\ln\left(\frac{1}{\varepsilon}\right)\right)$ | $\left\\|\vec{p}-\vec{p}_*\right\\|_\infty$ |
| MCMC | SG | $O\left(\frac{\ln n\ln(n/\sigma)}{\alpha\varepsilon^2}\right)$ | $\left\\|\vec{p}-\vec{p}_*\right\\|_2$ |

Поясним основные сокращения:

- *G-условие* – наличие такой web-страницы (например, страницы, отвечающей самой поисковой системе **G**oogle), на которую можно перейти из любой другой web-страницы с вероятностью не меньшей, чем $\alpha > cn^{-1}$, не ограничивая общности, будем считать, что вершина, отвечающая этой web-странице, имеет номер 1 (см. алгоритм MCMC);
- *S-условие* – из каждой web-страницы выходит в среднем не более $s \ll n$ ссылок на другие, то есть имеет место разреженность матрицы $P$ (**S**parsity); если из каждой web-страницы одновременно выходит и входит не более $s \ll n$ ссылок, то будем говорить об $\bar{S}$ условии;

  **Замечание 4.** В ряде случаев (например, для метода из [15]) можно релаксировать *S*-условие: “выходит в среднем”.
- *SG-условие* – спектральная щель $\alpha$ (**S**pectral **G**ap) матрицы $P$ удовлетворяет условию $\alpha \gg n^{-2}$, где $\alpha$ – расстояние между максимальным по величине модуля собственным значением (числом Фробениуса–Перрона) матрицы $P$ (равным 1) и модулем следую-щего (по величине модуля) собственного значения. Если выполняется G-условие, то выполняется и SG-условие с $\alpha$ не меньшим, чем в G-условии [2], [16].

Отметим, что приведенная таблица немного огрубляет результаты процитированных работ, в частности, работ [9], [10]. Это сделано для большей наглядности. Алгоритм из [8] на практике работает не очень быстро из-за большой константы в $\mathrm{O}(\cdot)$. Отчасти похожая ситуация и с алгоритмами из [11], [12], они также работают не так быстро, как можно было бы ожидать. Связано это также и с тем, что в стандартных пакетах (использовался MatLab) довольно не эффективно реализована возможность работы со случайностью в огромных размерностях. Алгоритмы из [14], [15] и MCMC работают в точности по приведенным оценкам. Заметно лучше приведенной оценки на практике работает алгоритм из [16], причем речь идет не о константе в $\mathrm{O}(\cdot)$. Тем не менее, условия, при которых этот алгоритм конкурентно способен довольно обременительные, и даже при этих условиях он, как правило, доминируем методом из [14]. Практический анализ всех приведенных в таблице (и не только) алгоритмов показал эффективность использования методов из [12], [15], когда мы не можем гарантировать то, что спектральная щель не мала. Если же мы можем это обеспечить, то неплохо работает метод MCMC при условии, что блок рандомизации пишется самостоятельно, т.е. не используется готовые способы генерирования дискретных случайных величин.

Обращает на себя внимание то, что у разных алгоритмов разные “цели”. В связи с этим полезно отметить, что “типично” (см. [17]):

$$\|\cdot\|_1 \sim \sqrt{n}\|\cdot\|_2,\ \|\cdot\|_2 \sim \sqrt{n}\|\cdot\|_\infty,$$

причем это соответствует векторам с одинаковыми по порядку компонентами. В нашем случае это как раз не типично (имеет место степенной закон убывания компонент (см. [18])), поэтому можно ожидать, что приведенные оценки будут с более слабой, чем $\sqrt{n}$ зависимостью, то есть при переходе от одной нормы к другой фактор $\sqrt{n}$ будет заменен чем-то более близким к $\mathrm{O}(1)$. Впрочем, мы не располагаем аккуратным обоснованием этого наблюдения.

По классификации главы 6 учебника [19] методы из [8], [9], [10] – являются вариационными, т.е. в этих методах решение системы $\vec{p}^T = \vec{p}^T P$ сводится к решению задач выпуклой оптимизации, которые, в свою очередь, решаются различными вариантами метода градиентного спуска (рандомизированный зеркальный спуск, метод Поляка–Шора, рандомизированный покомпонентный спуск).

Методы из [11], [12] – являются вариационными, но с игровым аспектом (оптимизационная задача понимается как задача поиска седловой точки – равновесия Нэша в матричной игре). Действительно, рассматриваемую нами задачу можно переписать (с помощью теории Фробениуса–Перрона [20]) как задачу:

$$f(\vec{p}) = \max_{\vec{u} \in S_n(1)} \langle \vec{u}, A\vec{p} \rangle \to \min_{\vec{p} \in S_n(1)},$$

где

$$A = P^T - I,\ A = \|a_{ik}\|,\ S_n(1) = \left\{ \vec{p} \geq \vec{0} :\ \sum_{k=1}^{n} p_k = 1 \right\}.$$

Отметим, что $0 \leq f(\vec{p}) \leq \|A\vec{p}\|_\infty$, $\vec{p} \in S_n(1)$, и на векторе PageRank (и только на нем) $f(\vec{p}) = 0$. Фактически, этот целевой функционал $f(\vec{p})$ оказывается очень близким к $\|P^T \vec{p} - \vec{p}\|_\infty$. В основе обоих игровых методов лежат рандомизированные варианты метода зеркального спуска поиска равновесия в матричных играх.

Методы из [8], [11], [12] могут быть распараллелены на $\sim \log_2(1/\sigma)$ процессорах. Если параллельно запустить $\sim \log_2(1/\sigma)$ независимых траекторий любого из этих методов, то, используя точные оценки числа итераций, гарантирующих заданную точность и доверительный уровень, который в этой параллельной схеме выбираем $1/2$, останавливаем траектории после выполнения предписанного числа итераций и проверяем выполнение критерия малости целевого функционала. С вероятностью $1-\sigma$ хотя бы одна из траекторий выдаст требуемое по точности решение. К сожалению, проверка выполнения критерия малости целевого функционала требует осуществления умножения $P^T \vec{p}$, что может быть более затратным, чем получить кандидата на решение. Поэтому, такое распараллеливание не всегда осмыслено. Кроме того, отмеченная проблема проверки выполнения критерия малости целевого функционала приводит к необходимости точного оценивания числа итераций, гарантирующих заданную точность и доверительный уровень (см., например, оценку на число итераций в теореме п. 5). Отметим, что метод из [12] без учета разреженности допускает распараллеливание на $n/\ln n$ процессорах (см. [12]).

Остальные методы из таблицы можно интерпретировать как вариации метода простых итераций из [19], [21]. При этом сразу же возникает вопрос: почему бы не попытаться решать систему $\vec{p}^T = \vec{p}^T P$ каким-то уже известным способом, например, методом из [19]? Или, скажем, таким (см. [22]): $\vec{p}^T = \vec{p}^T(0) P^\infty$, где $P^\infty$ – то, что выдает следующий итерационный процесс, требующий конечного числа итераций:

$$P(0) = I,\ P(k) = I - k\frac{P(k-1)(I-P)}{\operatorname{tr}\left(P(k-1)(I-P)\right)},\ k \in \mathbb{N}.$$

Итерации заканчиваются, когда первый раз выполнится условие $\operatorname{tr}\left(P(m)(I-P)\right) = 0$, где $\operatorname{tr}(B)$ – след матрицы $B$, т.е. сумма диагональных элементов. При этом $P^\infty = P(m)$.

Постараемся здесь ответить на этот вопрос. Прежде всего, многие итерационные методы требуют, чтобы матрица $P$ была симметричной и положительно определенной, что, как правило, место не имеет. Но даже если все необходимые условия будут выполнены и матрица при этом будет еще и разрежена, то будет наблюдаться такая же сходимость, как и в методе из [14], только вместо константы $\alpha$ будет фигурировать некий её аналог (из SG-условия), который в общем случае намного сложнее оценить. Таким образом, вместо того, чтобы приводить в таблице линейку классических итерационных методов с их оценками, мы ограничились тем, что привели наиболее приспособленные из этой линейки методы (условия применимости которых хорошо интерпретируемы) для решения рассматриваемого (довольно узкого) класса задач. Метод из [22], [23] мы не можем использовать, потому что $n \gg 10^4$ (т.е. $n^3 \gg 10^{12}$ – такое количество арифметических операций на одном персональном компьютере может быть выполнено за час).

В связи с упоминанием методов простой итерации, можно обратить внимание на то, что в описанных методах возможны проблемы накапливания ошибок округления (конечности длины мантиссы), возникающие при (см., например, § 4 главы 6 [19] или § 1 главы 5 [21]) условии $\|P\| > 1$ (при этом спектр матрицы может лежать в единичном круге, и метод простой итерации теоретически (без ошибок округления) должен сходиться со скоростью геометрической прогрессии). Однако в рассматриваемом нами случае в естественной норме, подчиненной $l_1^n$, $\|P\| = 1$, и таких проблем не возникает.

Важно также обратить внимание на то, что в ряде случаев, рассмотренных в таблице, целью является получение такого вектора распределения вероятностей $\vec{p}$, который давал бы малую невязку “по функции”. Но из того, что $\left\|P^T\vec{p}-\vec{p}\right\|$ малó не следует, что будет малó $\left\|\vec{p}-\vec{p}_*\right\|$.

**Замечание 5.** Напомним, что вектор $\vec{p}_*$ – решение уравнения: $P^T\vec{p}=\vec{p}$. Мы считаем это решение единственным в классе распределений вероятностей.

Более того, вполне естественно ожидать обратное, что $\left\|\vec{p}-\vec{p}_*\right\|$ окажется в $\alpha^{-1}$ раз больше (см. [19], [24]). Это оценка сверху, но, по-видимому, с большой вероятностью (если выбирать точку старта равновероятно) она, в действительности, превратится практически в точное соотношение (см. [24], [25]). Таким образом, добиться малости $\left\|P^T\vec{p}-\vec{p}\right\|$ – совсем не значит полностью решить задачу поиска вектора PageRank. Это обстоятельство, а также тот факт, что решение ищется на сильно ограниченном множестве (единичном симплексе) и дополнительно делаются всякие упрощающие предположения, отчасти объясняют, почему сложность описанных алгоритмов оказывается столь небольшой. Ведь, например, из того, что $\left\|\vec{p}-\vec{p}_*\right\|_\infty \le \varepsilon$ в случае, когда $\varepsilon \gg n^{-1}$ и истинное распределение $\vec{p}_*$ близко к равномерному, совсем не ясно, что именно выдает алгоритм с точностью $\varepsilon$ по функции, и как это можно использовать для ранжирования web-страниц. К счастью, такие ситуации, когда истинное распределение $\vec{p}_*$ близко к равномерному – нетипичны на практике (см. [27]), особенно в случае выполнения G-условия (SG-условия). Как правило, выделяются компоненты вектора $\vec{p}_*$, которые достаточно велики. А поскольку содержательно задача формулируется как ранжирование web-страниц, то по-сути, речь идет о восстановлении нескольких первых по величине компонент вектора $\vec{p}_*$. Другими словами, точно восстанавливать малые компоненты вектора $\vec{p}_*$ не требуется, если мы знаем, что они достаточно малы и есть достаточное количество не малых компонент. Скажем вектор, выдаваемый алгоритмом [16], имеет отличными от нуля не более чем $3\varepsilon^{-1}$ компонент. Отмеченные обстоятельства объясняют, почему даже при $\varepsilon \gg n^{-1}$ может быть полезен вектор $\vec{p}$, доставляющий оценку $\left\|P^T\vec{p}-\vec{p}\right\|_\infty \le \varepsilon$ – в одной из самых “плохих” норм: $l_\infty^n$.

## 3. Метод Markov Chain Monte Carlo

Идея решать линейные уравнения с помощью MCMC столь же стара, сколь и обычный метод Монте-Карло из [26]. Однако мы имеем дело не с линейным уравнением общего вида, а с уравнением со стохастической матрицей $P$, причем специальным образом заполненной, эти обстоятельства позволят нам более экономно организовать случайное блуждание по графу, соответствующему матрице $P$.

Прежде чем излагать алгоритм, приведем некоторые вспомогательные факты.

**Алгоритм Кнута–Яо (см. [26], [28]).** С помощью бросаний симметричной монетки требуется сгенерировать распределение заданной дискретной случайной величины (с.в.), принимающей конечное число значений. Предположим, что нам нужно сгенерировать распределение с.в., принимающей три значения 1, 2, 3 с равными вероятностями 1/3. Действуем таким образом. Два раза кидаем монетку: если выпало 00, то считаем что выпало значение 1, если 01, то 2, если 11, то 3. Если 10, то еще два раза кидаем монетку и повторяем рассуждения. Можно показать, что описанную выше схему можно так обобщить, чтобы сгенерировать распределение дискретной с.в., принимающей, вообще говоря, с разными вероятностями $n$ различных значений, в среднем с помощью не более чем $\log_2(n-1)+2$ подбрасываний симметричной монетки. Если эти вероятности одинаковы, то процедура "приготовления" такого алгоритма генерирования также имеет логарифмическую сложность по $n$.

**Алгоритм Markov Chain Monte Carlo (MCMC) (см. [7], [29]–[37]).** Чтобы построить однородный дискретный марковский процесс с конечным числом состояний, имеющий наперед заданную инвариантную (стационарную) меру $\pi$, переходные вероятности ищутся в следующем виде: $p_{ij} = p_{ij}^0 b_{ij}$, $i \neq j$; $p_{ii} = 1 - \sum_{j:\ j\neq i} p_{ij}$, где $p_{ij}^0$ – некоторая "затравочная" матрица, которую будем далее предполагать симметричной. Легко проверить, что матрица $p_{ij}$ имеет инвариантную (стационарную) меру $\pi$:

$$\frac{b_{ij}}{b_{ji}} = \frac{\pi_j p_{ji}^0}{\pi_i p_{ij}^0} = \frac{\pi_j}{\pi_i}, \ p_{ij}^0 > 0.$$

Чтобы найти $b_{ij}$ достаточно найти функцию $F:\mathbb{R}_+ \to [0,1]$ такую, что

$$\frac{F(z)}{F(1/z)} = z$$

и положить

$$b_{ij} = F\left(\frac{\pi_j p_{ji}^0}{\pi_i p_{ij}^0}\right) = F\left(\frac{\pi_j}{\pi_i}\right).$$

Пожалуй, самый известный пример такой функции $\tilde{F}(z) = \min\{z,1\}$ – алгоритм (Хастингса–)Метрополиса [7]. Заметим, что для любой такой функции $F(z)$ имеем $F(z) \le \tilde{F}(z)$. Другой пример дает функция $F(z) = z/(1+z)$. Заметим также, что $p_{ij}^0$ обычно выбирается равным $p_{ij}^0 = 1/M_i$, где $M_i$ число “соседних” состояний у $i$, или

$$p_{ij}^0 = 1/(2M),\ i \ne j;\ p_{ii}^0 = 1/2,\ i \ne j.$$

При больших значениях времени $t$, согласно эргодической теореме, имеем, что распределение вероятностей близко к стационарному $\pi$. Действительно, при описанных выше условиях имеет место условие детального баланса (марковские цепи, для которых это условие выполняется, иногда называют обратимыми):

$$\pi_i p_{ij} = \pi_j p_{ji},\ i, j = 1,...,n,$$

из которого сразу следует инвариантность меры $\pi$, т.е.

$$\sum_i \pi_i p_{ij} = \pi_j \sum_i p_{ji} = \pi_j,\ j = 1,...,n.$$

Основное применение замеченного факта состоит в наблюдении, что время выхода марковского процесса на стационарную меру (mixing time (см. [38])) во многих случаях оказывается удивительно малым.

**Замечание 6.** Более того, задача поиска такого симметричного случайного блуждания на графе (с равномерной инвариантной мерой в виду симметричности) заданной структуры, которое имеет “наименьшее” mixing time (другими словами, наибольшую спектральную щель), сводится к задаче полуопределенного программирования, которая, как известно, полиномиально (от числа вершин этого графа) разрешима [39]

При том что выполнение одного шага по времени случайного блуждания по графу, отвечающему рассматриваемой марковской цепи, как следует из алгоритма Кнута–Яо, также может быть быстро сделано. Таким образом довольно часто можно получать эффективный способ генерирования распределения дискретной случайной величины с распределением вероятностей $\pi$ за время, полиномиальное от логарифма числа компонент вектора $\pi$.

Для оценки mixing time нужно оценить спектральную щель стохастической матрицы, переходных вероятностей, задающей исследуемую марковскую динамику, то есть нужно оценить расстояние от максимального собственного значения этой матрицы равного единицы (теорема Фробениуса–Перрона) до следующего по величине модуля. Именно это число определяет основание геометрической прогрессии, мажорирующей исследуемую последовательность норм разностей расстояний (по вариации) между распределением в данный момент времени и стационарным (финальным) распределением. Для оценки спектральной щели разработано довольно много методов, из которых мы упомянем лишь некоторые (см. [30]–[37]): неравенство Пуанкаре (канонический путь), изопериметрическое неравенство Чигера (проводимость), с помощью техники каплинга [40] (получаются простые, но, как правило, довольно грубые оценки), с помощью каплинга Мертона [41], с помощью дискретной кривизны Риччи и теорем о концентрации меры (Мильмана–Громова [37], [42]). Приведем некоторые примеры применения MCMC [7]: Тасование $n$ карт, разбиением приблизительно на две равные кучи и перемешиванием этих куч (mixing time $\sim \log_2 n$); Hit and Run (mixing time $\sim n^3$); Модель Изинга – $n$ спинов на отрезке, стационарное распределение = распределение Гиббса, Глауберова динамика (mixing time $\sim n^{2\log_2 e/T}$, $0 < T \ll 1$); Проблема поиска кратчайших гамильтоновых путей; Имитация отжига для решения задач комбинаторной оптимизации, MCMC для решения задач перечислительной комбинаторики.

**Замечание 7.** На примере тасования карт контраст проявляется, пожалуй, наиболее ярко. Скажем для колоды из 52 карт пространство состояний марковской цепи будет иметь мощность 52! (если сложить времена жизней в наносекундах каждого человека, когда либо жившего на Земле, то это число на много порядков меньше 52!). В то время как такое тасование: взять сверху колоды карту, и случайно поместить ее во внутрь колоды, отвечающее определенному случайному блужданию, с очень хорошей точностью выйдет на равномерную меру, отвечающую перемешанной колоде, через каких-то 200–300 шагов. Если брать тасование разбиением на кучки, то и того меньше – за 8–10 шагов [45].

Продемонстрируем сказанное выше двумя примерами (с помощью которых, например, получаются оценки mixing time работы [43]), которые нам пригодятся в дальнейшем.

**Пример 1 (подход Чигера [31], [34], [36]).** Пусть (напомним, что под $\pi(\cdot)$ понимается инвариантная мера, а под $P = \left\| p_{ij} \right\|_{i,j=1}^{n,n}$ – матрица переходных вероятностей марковской цепи)

$$h(G) = \min_{S\subseteq V_G:\pi(S)\le 1/2} P\left(S \to \bar{S} \middle| S\right) = \min_{S\subseteq V_G:\pi(S)\le 1/2} \frac{\sum\limits_{(i,j)\in E_G:i\in S,\, j\in\bar{S}} \pi(i)\, p_{ij}}{\sum\limits_{i\in S}\pi(i)}, \quad \text{(константа Чигера)}$$

$$T(i,\varepsilon) = \Theta\left(h(G)^{-2}\left(\ln\left(\pi(i)^{-1}\right) + \ln\left(\varepsilon^{-1}\right)\right)\right). \quad \text{(Mixing time)}$$

Тогда

$$\forall \; i = 1,\ldots,n, \; t \ge T(i,\varepsilon) \to \left\| P^t(i,\cdot) - \pi(\cdot) \right\|_1 = \sum_{j=1}^{n} \left| P^t(i,j) - \pi(j) \right| \le \varepsilon,$$

где $P^t(i,j)$ – условная вероятность того, что стартуя из состояния $i$ через $t$ шагов, марковский процесс окажется в состоянии $j$. Отметим, что от вектора PageRank $\pi$ мы вправе ожидать степенной закон убывания компонент, отсортированных по возрастанию (см. [27]), поэтому $\max\limits_{i=1,\ldots,n} \ln\left(\pi(i)^{-1}\right) = O(\ln n)$.

**Пример 2 (coarse Ricci curvature [37]).** Введем расстояние Монжа–Канторовича между двумя (дискретными) распределениями вероятностей $\mu$ и $\nu$:

$$W_1(\mu,\nu) = \min_{\substack{\xi \ge 0:\, \sum_j \xi(i,j) = \mu(i) \\ \sum_i \xi(i,j) = \nu(j)}} \sum_{i,j} d(i,j)\,\xi(i,j),$$

где каждой паре вершин поставлено в соответствие неотрицательное число $d(i,j)$ (со свойствами расстояния (метрики)). Говорят, что $\kappa$ – дискретная кривизна Риччи, если

$$\exists \; t_0 > 0: \; \forall \; i,j = 1,\ldots,n$$
$$W_1\left(P^{t_0}(i,\cdot), P^{t_0}(j,\cdot)\right) \le (1-\kappa)\, d(i,j).$$

Пусть существует такое $\kappa > 0$, тогда

$$W_1\left(P^t(i,\cdot), \pi(\cdot)\right) \le \left[\kappa^{-1} \sum_{j=1}^{n} d(i,j) p_{ij}\right] \cdot (1-\kappa)^{t/t_0},$$

$$E\left[\left\|\frac{1}{T}\sum_{t=T_0}^{T+T_0}\vec{x}_t-\vec{\pi}\right\|_2^2\right]=\left(\text{Bias}\right)^2+\text{Var}=\left(O\left(n\frac{\left(1-\kappa\right)^{T_0}}{\kappa T}\right)\right)^2+O\left(\frac{1}{\kappa T}\right),$$

где у случайного вектора $\vec{x}_t$ все компоненты нулевые, кроме единственной компоненты, отвечающей вершине, в которой находился марковский случайный процесс на шаге $t$. Считая, что

$$\text{Var}\gg\left(\text{Bias}\right)^2$$

имеем (этот результат в явном виде не содержится в [37], но к нему можно прийти, используя некоторые идеи работы [41])

$$\exists\ c>0:\ P\left(\left\|\frac{1}{T}\sum_{t=T_0}^{T+T_0}\vec{x}_t-\vec{\pi}\right\|_2>c\sqrt{\frac{\ln n+\Omega}{\kappa T}}\right)\le\exp\left(-\Omega\right).$$

Результат вида

$$\frac{1}{T}\sum_{t=T_0}^{T+T_0}\vec{x}_t\xrightarrow[T\to\infty]{}\vec{\pi}$$

– есть эргодическая теорема для марковских процессов (при этом выше была приведена довольно тонкая оценка того, какая скорость сходимости), вполне ожидаем (см. [15]) и соответствует классическим вариантам эргодических теорем для динамических системах (Биркгофа–Хинчина, фон Неймана). Правда, в отличие от динамических систем здесь удается оценить скорость сходимости.

Отметим связь описанного в примере 2 подхода с результатами о сжимаемости (в пространстве всевозможных лучей с центром в начале координат неотрицательного ортанта) в метрике Биркгофа положительных линейных операторов (см. [44]), в частности, заданных стохастической матрицей $P$, некоторая степень которой имеет все элементы положительными (это равносильно неразложимости и непериодичности марковской цепи (см. [45])). Кстати сказать, такое понимание эргодической теоремы для марковских цепей позволяет интерпретировать её как теорему о сжимающих отображениях, что выглядит несколько необычно в контексте сопоставления этой теоремы с эргодическими теоремами для динамических систем.

# 4. Алгоритм MCMC

Теперь можно приступать к изложению нужной нам версии алгоритма MCMC.

**Шаг 1. Инициализация:** $\vec{X}=\vec{0}$, вершина $=1$, $t=0$.

**Шаг 2. Счетчик итераций:** $t:=t+1$.

**Шаг 3. Модификация** $\vec{X}$ **:** если $t\ge T^{0}_{\varepsilon,\alpha,n}$, то $X_k:=X_k+1$, где $k$ – номер текущей вершины.

**Шаг 4. Модификация вершины:** случайно “переходим” из текущей вершины в одну из “соседних” согласно матрице $P$.

**Шаг 5. Остановка:** если $t<T_{\varepsilon,\sigma,\alpha,n}$ перейти на шаг 2, иначе на шаг 6.

**Шаг 6. Ответ:** $\vec{p}=\vec{X}\Big/\left(T_{\varepsilon,\sigma,\alpha,n}-T^{0}_{\varepsilon,\alpha,n}\right)$.

Самым затратным шагом из первых пяти шагов является шаг 4, но даже этот шаг при самом неблагоприятном раскладе (из текущей вершины выходит порядка $n$ ребер) можно осуществить за $O(\ln n)$ операций (см., например, алгоритм Кнута–Яо – на самом деле все это можно сделать разными способами, причем сложность $O(\ln n)$ можно понимать не в среднем, а обычном образом (см. [26]); алгоритм Кнута–Яо был приведен лишь как иллюстративный пример). Из примеров 1 и 2 (мы заменяем оценки спектральной щели в этих примерах на саму спектральную щель $\alpha$, что можно делать в таком контексте для обратимых марковских цепей (см. [41]), и нуждается в некоторых оговорках в общем случае), получаем, что при

$$T^{0}_{\varepsilon,\alpha,n}=O\left(\frac{1}{\alpha}\ln\left(\frac{n}{\varepsilon}\right)\right)$$

после

$$T_{\varepsilon,\sigma,\alpha,n}=O\left(\frac{\ln\left(n/\sigma\right)}{\alpha\varepsilon^{2}}\right)$$

итераций с вероятностью не меньшей $1-\sigma$ алгоритм MCMC выдаст $\varepsilon$-оптимальное (по функции) решение исходной задачи. При этом алгоритм MCMC затрачивает в общей сложности

$$\mathrm{O}\left(n+\frac{\ln n\ln(n/\sigma)}{\alpha\varepsilon^2}\right)$$

элементарных арифметических операций (типа умножения двух чисел с плавающей точкой) – слагаемое $\mathrm{O}(n)$ “отражает” стоимость шага 6.

Несложно предложить способ эффективного распараллеливания такого алгоритма. Для этого выпускается не одна траектория, а

$$N_{\varepsilon,\sigma}=\mathrm{O}\left(\frac{\ln(\sigma^{-1})}{\varepsilon^2}\right),$$

независимо блуждающих, и стартующих с вершин, выбираемых случайно и независимо. За каждой траекторией следим время $T^0_{\varepsilon/2,\alpha,n}$. Потом усредняем (см. ниже) по всем этим траекториям, получаем ответ (с возможностью организации параллельных вычислений), но при этом придется затратить чуть больше операций – вместо множителя $\ln(n/\sigma)$ появится $\ln(n/\varepsilon)\ln(\sigma^{-1})$, что не так уж и плохо.

**Замечание 8.** Более того, это обстоятельство на практике можно частично нивелировать, например, таким образом. Сначала выпускается одна траектория. Блуждающая частица через случайные моменты времени “рождает потомков”, которые также начинают независимо блуждать, стартуя с “места рождения”, причем рожденные частицы также склонны к “спонтанному делению”. Здесь важно правильно подобрать интенсивность такого деления (размножения).

Действительно, мы получаем набор из $N_{\varepsilon,\sigma}$ независимых одинаково распределенных случайных векторов $\vec{x}^k_{T^0_{\varepsilon,\alpha,n}}$, $k=1,\ldots,N_{\varepsilon,\sigma}$, где $\vec{x}^k_{T^0_{\varepsilon,\alpha,n}}$ – вектор, все компоненты которого равны нулю кроме одной, равной 1; эта компонента соответствует состоянию, в котором находится $k$-е блуждание на шаге $T^0_{\varepsilon/2,\alpha,n}$. Векторы одинаково распределены, причем (напомним, что $\vec{p}_*=\vec{\pi}$ – инвариантная мера)

$$\left\|E\left[\vec{x}^k_{T^0_{\varepsilon/2,\alpha,n}}-\vec{\pi}\right]\right\|_1\le\varepsilon/2.$$

Применяем далее (этот способ был предложен совместно с Е.Ю. Клочковым [46]) к

$$\frac{1}{N_{\varepsilon,\sigma}}\sum_{k=1}^{N_{\varepsilon,\sigma}}\vec{x}^{k}_{T^0_{\varepsilon/2,\alpha,n}}$$

неравенство типа Хефдинга в гильбертовом пространстве $l_2^n$ (см. [49]), получаем

$$P\left(\left\|\frac{1}{N_{\varepsilon,\sigma}}\sum_{k=1}^{N_{\varepsilon,\sigma}}\vec{x}^{k}_{T^0_{\varepsilon/2,\alpha,n}}-\vec{\pi}\right\|_2\geq\varepsilon\right)=P\left(\left\|\sum_{k=1}^{N_{\varepsilon,\sigma}}\left(\vec{x}^{k}_{T^0_{\varepsilon/2,\alpha,n}}-\vec{\pi}\right)\right\|_2\geq\varepsilon N_{\varepsilon,\sigma}\right)\leq$$

$$\leq\exp\left(-\frac{1}{4N_{\varepsilon,\sigma}}\left(\varepsilon N_{\varepsilon,\sigma}-E\left(\left\|\sum_{k=1}^{N_{\varepsilon,\sigma}}\left(\vec{x}^{k}_{T^0_{\varepsilon/2,\alpha,n}}-\vec{\pi}\right)\right\|_2\right)\right)^2\right),$$

где

$$E\left(\left\|\sum_{k=1}^{N_{\varepsilon,\sigma}}\left(\vec{x}^{k}_{T^0_{\varepsilon/2,\alpha,n}}-\vec{\pi}\right)\right\|_2\right)\leq\sqrt{E\left(\left\|\sum_{k=1}^{N_{\varepsilon,\sigma}}\left(\vec{x}^{k}_{T^0_{\varepsilon/2,\alpha,n}}-\vec{\pi}\right)\right\|_2^2\right)}\leq\sqrt{2\sqrt{2}N_{\varepsilon,\sigma}+\left(\varepsilon^2/4\right)N_{\varepsilon,\sigma}^2}.$$

Подберем $N_{\varepsilon,\sigma}$ так, чтобы выполнялось неравенство

$$P\left(\left\|\frac{1}{N_{\varepsilon,\sigma}}\sum_{k=1}^{N_{\varepsilon,\sigma}}\vec{x}^{k}_{T^0_{\varepsilon/2,\alpha,n}}-\vec{\pi}\right\|_2\geq\varepsilon\right)\leq\sigma.$$

Для этого достаточно, чтобы выполнялось неравенство

$$\varepsilon N_{\varepsilon,\sigma}-\sqrt{2\sqrt{2}N_{\varepsilon,\sigma}+\left(\varepsilon^2/4\right)N_{\varepsilon,\sigma}^2}\geq\sqrt{4N_{\varepsilon,\sigma}\ln\left(\sigma^{-1}\right)},$$

где

$$N_{\varepsilon,\sigma}=\frac{4+6\ln\left(\sigma^{-1}\right)}{\varepsilon^2}=O\left(\frac{\ln\left(\sigma^{-1}\right)}{\varepsilon^2}\right).$$

Отметим, что с точностью до мультипликативной константы эта оценка не может быть улучшена. Это следует из неравенства Чебышёва

$$P\left(X\geq EX-\varepsilon\right)\geq 1-\frac{\mathrm{Var}\left(X\right)}{\varepsilon^2},$$

при

$$X = \left\| \frac{1}{N_{\varepsilon,\sigma}} \sum_{k=1}^{N_{\varepsilon,\sigma}} \vec{x}^{k}_{T^0_{\varepsilon/2,\alpha,n}} - \vec{\pi} \right\|_2, \quad \vec{\pi} = \left(n^{-1},\ldots,n^{-1}\right)^T, \quad N_{\varepsilon,\sigma} \gg n.$$

**Замечание 9.** Стоит отметить, что при экспоненциальном убывании компонент отранжированного вектора $\vec{\pi}$ можно заменить полученные далее оценки в 2-норме на аналогичные оценки (с другими константами и дополнительным логарифмическим штрафом), но в 1-норме. Понять это можно из того, что фактически речь идет о задаче восстановления параметров мультиномиального распределения. Точнее, речь идет о влиянии размерности пространства параметров на оценки скорости сходимости в теореме Фишера о методе наибольшего правдоподобия в неасимптотическом варианте [47] (выборочные средние как раз и будут оценками, полученными согласно этому методу). Для данного примера мультиномиальное распределение порождает (вот в этом месте, к сожалению, приходится использовать формулу Стирлинга, то есть требуются некоторые асимптотические оговорки) в показателе экспоненты расстояние Кульбака–Лейблера между выборочным средним и истинным распределением [48]. В свою очередь это расстояние оценивается согласно неравенству Пинскера квадратом 1-нормы [49].

Немного удивляет крайне слабая зависимость общей сложности от размера матрицы $P$, особенно, если учесть, что никаких предположений о разреженности $P$ не делалось. “Подвох” здесь в том, что, в действительности, если речь не идет о каких-то специальных графах (например, экспандерах [50]), на которых рассматривается случайное блуждание, то условие, что $\alpha$ равномерно по $n$ отделимо от нуля – противоестественно. Так, в работе [43] приводится следующая, по-видимому, не улучшаемая оценка $\alpha \sim n^{-1}$ для класса часто встречающихся на практике марковских процессов, которые возникают при описании различных макросистем, “живущих” в приближении среднего поля (в динамику заложено равноправие взаимодействующих агентов – закон действующих масс Гульдберга–Вааге). Также можно заметить, что выигрыш в оценке сложности в зависимости от $n$ происходит за счет довольно плохой зависимости от $\varepsilon$. Вообще ожидать зависимости сложности от $\varepsilon$ лучшей, чем $\varepsilon^{-2}$ в рандомизированных алгоритмах не приходится (см. [51]), поэтому рандомизация осмыслена, как правило, только при $\varepsilon \gg n^{-1}$. Это становится особенно ясно, если сравнить полученную оценку сложности с оценкой сложности, скажем, алгоритма из [14]. То есть, чтобы MCMC имело смысл использовать нужно, чтобы $\varepsilon \gg n^{-1}$.

Основным недостатком метода MCMC является отсутствие точного знания о $T_{\varepsilon,\sigma,\alpha}$ и $T^0_{\varepsilon,\alpha,n}$, даже если известен размер спектральной щели $\alpha$. Более того, эффективность алгоритма напрямую завязана на оценку спектральной щели $\alpha$, которая, как правило, априорно не известна. Тем не мене, проблема оценки $\alpha$ может быть решена за $Cn$ арифметических операций (к сожалению, с довольно большой константой $C$), например, с помощью $\delta^2$-процесса практической оценки спектральной щели (см., например, [19]). Другой способ оценки $\alpha$ и $T^0_{\varepsilon,\alpha,n}$ имеется в [52]. Проблема определения $T_{\varepsilon,\sigma,\alpha}$ решается с помощью контроля разности $\|\vec{p}_{t+\tau} - \vec{p}_t\|_2$. Относительно $T^0_{\varepsilon,\alpha,n}$ можно действовать так: изначально запускать алгоритм с $T^0_{\varepsilon,\alpha,n} = 0$, а затем скорректировать полученный ответ, полагая, скажем, $T^0_{\varepsilon,\alpha,n} = T_{\varepsilon,\sigma,\alpha}/5$.

Отметим также, что изначально необходимо правильным образом разместить в памяти компьютера матрицу $P$. Если работать с этой матрицей обычным образом, то при случайном выборе алгоритмом MCMC на каждом шаге соседней вершины будет тратиться время порядка $n$, а не $\ln n$, как хотелось бы. Чтобы избежать этого, необходимо перед началом работы алгоритма представить граф в виде списка ссылок. Это займет время (память) порядка числа ребер, но сделать это нужно всего один раз.

## 5. Алгоритм Григориадиса–Хачияна

Рассматриваемую нами задачу перепишем (с помощью теории Фробениуса–Перрона [20]) как задачу поиска седловой точки (равновесия в антагонистической матричной игре):

$$f(\vec{p}) = \max_{\vec{u}\in S_n(1)} \langle \vec{u}, A\vec{p}\rangle \to \min_{\vec{p}\in S_n(1)}, \tag{1}$$

где $A = P^T - I$, $S_n(1) = \left\{\vec{p} \geq \vec{0}:\ \sum_{k=1}^{n} p_k = 1\right\}$. Отметим, что $f(\vec{p}) \geq 0$, $\vec{p} \in S_n(1)$, и на векторе PageRank (и только на нем) $f(\vec{p}) = 0$.

С точностью (по функции) $\varepsilon$ и с вероятностью не меньшей $1-\sigma$ равновесие можно найти с помощью стохастического зеркального спуска из работы [11] за количество операций

$$O\left(\frac{n\ln(n/\sigma)}{\varepsilon^2}\right).$$

Однако эти методы не достаточно полным образом учитывают специфику задачи в случае разреженной матрицы $P$. Для решения задачи (1) воспользуемся методом поиска равновесий в симметричных антагонистических матричных играх из [12] (отметим, что этот метод по сути является определенным образом рандомизированным методом зеркального спуска из [51]). Для этого, предварительно приведем задачу (1), следуя Данцигу [12], к симметричному виду (этого можно и не делать, если использовать рандомизированный онлайн метод зеркального спуска из теоремы 2 работы [53]):

$$\max_{\vec{u}\in S_{2n+1}(1)} \langle \vec{u}, \mathrm{A}\vec{x} \rangle \to \min_{\vec{x}:=(\vec{y},\vec{p}',u)\in S_{2n+1}(1)}, \quad \mathrm{A} := \begin{bmatrix} 0 & A & -\vec{e} \\ -A^T & 0 & \vec{e} \\ \vec{e}^T & -\vec{e}^T & 0 \end{bmatrix},$$

где $\vec{e} = (1,\ldots,1)^T$, $\mathrm{A} = \|a_{ik}\|$. Тогда, $f(\vec{p}) \le 2\varepsilon$, где $\vec{p} = \vec{p}'/(\vec{e}^T\vec{p}')$, причем $\vec{e}^T\vec{p}' \ge 1/2 - \varepsilon$, если $\mathrm{A}\vec{x} \le \varepsilon\vec{e}$.

**Шаг 1. Инициализация:** $\vec{X} = \vec{0}$, $\vec{p}^T = \dfrac{1}{2n+1}\underbrace{(1,\ldots,1)}_{2n+1}^T$, $t = 0$.

**Шаг 2. Счетчик итераций:** $t := t+1$.

**Шаг 3. Датчик случайных чисел:** выбираем $k \in \{1,\ldots,2n+1\}$ с вероятностью $p_k$.

**Шаг 4. Модификация** $\vec{X}$ **:** $X_k := X_k + 1$.

**Шаг 5. Модификация** $\vec{p}$ **:** $i = 1,\ldots,2n+1$ $p_i := p_i \exp\left(\dfrac{\varepsilon a_{ik}}{2}\right)$.

**Шаг 6. Остановка:** если $t < T_{\varepsilon,\sigma,n}$ перейти на шаг 2, иначе на шаг 7.

**Шаг 7. Ответ:** $\vec{x} = \vec{X}/t$.

**Замечание 10 (к шагу 5).** Далее будет пояснено, что достаточно задавать распределение вероятностей с точностью до нормирующего множителя.

Итак, считая матрицу $P$ разреженной с не более чем $s$ элементами в столбце и строке (см. [9]), получаем, что самыми “дорогими” шагами будет пересчет (перегенерирование) распределения вероятностей (шаги 3 и 5). Используя способ генерирования дискретной случайный величины с помощью сбалансированного двоичного дерева, мы получаем, что изменение веса вероятности одного исхода, по сути, равносильно процедуре изменение весов тех вершин дерева, путь через которые ведет к изменяемому значению “листа” дерева.

**Замечание 11.** Опишем точнее эту процедуру. У нас есть сбалансированное двоичное дерево высоты $O(\log_2 n)$ с $2n+1$ листом (дабы не вдаваться в технические детали, считаем что число $2n+1$ – есть степень двойки, понятно, что это не так, но можно взять, скажем, наименьшее натуральное $m$ такое, что $2^m > 2n+1$ и рассматривать дерево с $2^m$ листом) и с $O(n)$ общим числом вершин. Каждая вершина дерева (отличная от листа) имеет неотрицательный вес равный сумме весов двух ее потомков. Первоначальная процедура приготовления дерева, отвечающая одинаковым весам листьев, потребует $O(n)$ операций. Но сделать это придется всего один раз. Процедура генерирования дискретной случайной величины с распределением, с точностью до нормирующего множителя, соответствующим весам листьев может быть осуществлена с помощью случайного блуждания из корня дерева к одному из листьев. Отметим, что поскольку дерево двоичное, то прохождение каждой его вершины при случайном блуждании, из которой идут два ребра в вершины с весами $a>0$ и $b>0$, осуществляется путем подбрасывания монетки (“приготовленной” так, что вероятность пойти в вершину с весом $a$ – есть $a/(a+b)$). Понятно, что для осуществления этой процедуры нет необходимости в дополнительном условии нормировки: $a+b=1$. Если вес какого-то листа алгоритму необходимо поменять по ходу работы, то придется должным образом поменять дополнительно веса тех и только тех вершин, которые лежат на пути из корня к этому листу. Это необходимо делать, чтобы поддерживать свойство: каждая вершина дерева (отличная от листа) имеет вес равный сумме весов двух ее потомков.

Обратим при этом внимание, что нет необходимости заниматься перенормировкой распределения вероятностей (это бы стоило $O(n)$), то есть изменением весов вершин дерева, отличных от тех, путь через которые ведет к изменяемому значению листа дерева. Все это (сгенерировать с помощью дерева и обновить дерево) можно сделать за $O(\ln n)$ операций типа сравнения двух чисел, а в (“типичном”) случае, когда это нужно делать $\sim s$ раз, то за $O(s\ln n)$ операций (и лишь в случаях, когда $k=2n+1$, придется делать $\sim 2n$ операций). Далее будет показано, что после

$$T_{\varepsilon,\sigma,n} = \mathrm{O}\left(\frac{\ln(n/\sigma)}{\varepsilon^2}\right)$$

итераций с вероятностью, не меньшей $1-\sigma$ алгоритм выдаст $\varepsilon$-оптимальное (по функции) решение исходной задачи. При этом алгоритм затрачивает в общей сложности

$$O\left(n+\frac{s\ln n\ln(n/\sigma)}{\varepsilon^2}\right)$$

элементарных арифметических операций (типа умножения двух чисел с плавающей точкой), что в случае $s \ll n$ заметно лучше всех известных сейчас оптимизационных аналогов. Слагаемое $O(n)$ “отражает” те случаи, когда $k=2n+1$, а также “стоимость” заключительного шага 7. Таким образом, имеет место следующее утверждение, доказательство которого вынесено в приложение.

**Теорема.** *Алгоритм Григориадиса–Хачияна после*

$$T_{\varepsilon,\sigma,n} = 12\left(\ln(2n+1) + \ln\left(\sigma^{-1}\right)\right)\varepsilon^{-2}$$

*итераций с вероятностью, не меньшей* $1-\sigma$ *выдает такое* $\vec{p}$*, что* $0 \le f(\vec{p}) \le \varepsilon$*. При этом затрачивается в общей сложности*

$$O\left(n + \frac{s \ln n \ln(n/\sigma)}{\varepsilon^2}\right)$$

*операций (вида умножения двух чисел типа double).*

Ранее мы уже отмечали, что улучшить зависимость сложности от $\varepsilon$ не представляется возможным, существенно не ухудшая зависимость сложности от $n$. Можно даже сказать точнее, что для методов, в которых используется рандомизация с дисперсией стохастического градиента $D$, зависимость числа итераций $\sim D\varepsilon^{-2}\ln\left(\sigma^{-1}\right)$ типична и не улучшаема; тем не менее, стоит оговориться, что если от сильно выпуклого функционала берется математическое ожидание (это не наш случай, мы, напротив, избавлялись от этого, чтобы функционал был более информативен), то сходимость может быть улучшена до $\sim D\varepsilon^{-1}$ (см. [54], [55]). Например, если использовать (детерминированный) быстрый градиентный метод (см. [56], [57]), то можно получить зависимость сложности от $\varepsilon$ вида $\varepsilon^{-1}$, но при этом число операций увеличится не менее чем в $n$ раз. Поскольку для таких задач вполне естественным является соотношение $\varepsilon \gg n^{-1}$, то выгода от этого представляется сомнительной. Отметим также (см. [12]), что в классе детерминированных алгоритмов зависимость от $n$ не может быть лучше, чем $\sim n^2$ (речь идет не о разреженных матрицах $P$).

Если бы мы применили обычный метод зеркального спуска (см. [58], [59]) для поиска (безусловного) минимума негладкой выпуклой функции $f(\vec{p}) + \left(\sum_{k=1}^{n} p_k - 1\right)^2$, то используя теорему о субдиференциале максимума из [60], [61], и считая матрицу $A$ слабо заполненной, со средним показателем заполненности (по строкам и столбцам) $\chi \ll 1$ (в частности, число элементов матрицы отличных от нуля будет $\le \chi n^2$), получили бы оценку общей сложности в среднем

$$O\left(\chi n^2+\frac{\left(n+\chi^2 n^2\right)\ln n}{\varepsilon^2}\right).$$

В определенных ситуациях (например, при $\chi \sim n^{-1/2}$) эта оценка оказывается вполне конкурентноспособной.

Отметим связь описанного алгоритма с онлайн оптимизацией. Грубо говоря, такого типа алгоритмы асимптотически наиболее эффективны в задачах онлайн обучения на основе опыта экспертов даже при сопротивляющейся “Природе” (см. [51], [53], [54], [55], [59], [61], [62], [63], [64]).

Отметим также связь алгоритма Григориадиса–Хачияна с концепцией ограниченной рациональности в контексте Discrete choice theory (см. [53], [62], [65]).

Резюмируем полученные в этом пункте результаты. Описанный выше алгоритм Григориадиса–Хачияна фактически соответствуют методу зеркального спуска для (1) с кососимметричной матрицей $\mathrm{A}$. Для не кососимметричной матрицы см. [53].

Отличие рассмотренного зеркального спуска от зеркального спуска, например, работы [11] в том, что в нашей работе (также как и в [12] движение итерационного процесса для поиска седловой точки $\langle\vec{u},\mathrm{A}\vec{x}\rangle$ на произведении двух единичных симплексов осуществляется только по $\vec{x}$, в то время как обычный метод зеркального спуска (рассчитанный не на кососимметричную матрицу $\mathrm{A}$) осуществляет движение также и по $\vec{u}$. Кроме того, очень важно, как именно осуществлять рандомизацию. В этой работе, также как и в [12], рандомизация происходит при выборе компоненты вектора $\vec{x}$, по которой осуществляется покомпонентное движение, во всех других работах рандомизацию предлагается вводить на этапе вычисления $\mathrm{A}\vec{x}=E\left[\mathrm{A}^{\langle x\rangle}\right]$, где $x$ – дискретная с.в. с распределением $\vec{x}$, а $\mathrm{A}^{\langle x\rangle}$ есть $x$-столбец матрицы $\mathrm{A}$. Оба описанных способа рандомизации требуют одинакового по порядку числа шагов

$$O\left(\frac{\ln\left(n/\sigma\right)}{\varepsilon^2}\right)$$

для достижения точности (по функции) $\varepsilon$, но рандомизация Григориадиса–Хачияна позволяет заметно лучше учитывать разреженную специфику. Причина проста – на каждом шаге рандомизированных зеркальных спусков работ [8], [11] требуется обновлять вектор $\vec{x}$, прибавляя к нему вектор с $2n+1$ ненулевыми компонентами, в независимости от раз-реженности А, стало быть, тратить не менее $2n+1$ операций на один шаг. В то время как у алгоритма Григориадиса–Хачияна в разреженном случае число операций на шаг оказывается в типичной ситуации, как мы видели выше, $O(s\ln n)$, что может быть заметно меньше $n$ в случае $s \ll n$.

В перспективе планируется развивать идею о том, как организовывать случайный покомпонентный градиентный спуск для задач выпуклой оптимизации в пространствах огромной размерности, чтобы как можно сильнее учесть разреженную специфику задачи. Сейчас популярным являются другие подходы (см. [9], [10]). Однако нам представляется, что в определенных ситуациях описанный здесь способ рандомизации (Григориадиса–Хачияна) будет давать лучшие результаты. Также планируется исследовать вопросы о робастном оценивании вектора PageRank (см. [66], [67]) и развить некоторые идеи, связанные с распределенными вычислениями вектора PageRank (см. [68]).

## Приложение: Доказательство теоремы

Докажем теорему, сформулированную выше. Далее (для простоты будем вместо $2n+1$ писать $n$). Определим $r>0$ из соотношения: $\sigma \simeq n^{-r}$. Покажем, во многом следуя работе [12], что алгоритм Григориадиса–Хачияна после

$$T_{\varepsilon,\sigma,n} = 3\left(\ln n + \ln\left(\sigma^{-1}\right)\right)\varepsilon^{-2} = 3(1+r)\varepsilon^{-2}\ln n$$

итераций выдает такое $\vec{x}$, что с вероятностью не меньшей чем $1-\sigma$, имеет место неравенство

$$\mathrm{A}\vec{x} \le \varepsilon\vec{e}.$$

Сначала, следуя работам [12], [61], [62], положим

$$p_i(t) = P_i(t)\left(\sum_{j=1}^{n} P_j(t)\right)^{-1}, \quad P_i(t) = \exp\left(\varepsilon U_i(t)/2\right)$$

и

$$\Phi(t) = \sum_{i=1}^{n} P_i(t), \text{ где } \vec{U}(t) = \mathrm{A}\vec{X}(t).$$

Далее, аналогично [12], имеем

$$\Phi(t+1) = \sum_{i=1}^{n} P_i(t)\exp\left(\varepsilon a_{ik}/2\right) =$$

$$= \Phi(t)\sum_{i=1}^{n} p_i(t)\exp(\varepsilon a_{ik}/2),$$

$$E\left[\Phi(t+1)\middle|\vec{P}(t)\right] = \Phi(t)\sum_{i,k=1}^{n} p_i(t)p_k(t)\exp(\varepsilon a_{ik}/2)$$

и ($|a_{ik}| \le 1$)

$$\exp(\varepsilon a_{ik}/2) \le 1 + \varepsilon a_{ik}/2 + \varepsilon^2/6,$$

но поскольку

$$\sum_{i,k=1}^{n} p_i(t)\, p_k(t) = \left(\sum_{i=1}^{n} p_i(t)\right)^2 = 1,\ \sum_{i,k=1}^{n} p_i(t)\, p_k(t)\frac{\varepsilon^2}{6} = \frac{\varepsilon^2}{6},$$

$$\sum_{i,k=1}^{n} p_i(t)\, p_k(t) a_{ik} = \left\langle \vec{p}(t), \mathrm{A}\vec{p}(t)\right\rangle = 0,$$

то

$$E\left[\Phi(t+1)\middle|\vec{P}(t)\right] \le \Phi(t)\left(1+\varepsilon^2/6\right),$$

$$E\left[\Phi(t+1)\right] \le E\left[\Phi(t)\right]\left(1+\varepsilon^2/6\right).$$

Используя это неравенство и то, что

$$E\left[\Phi(0)\right] = \Phi(0) = n,$$

имеем

$$E\left[\Phi(t)\right] \le n\left(1+\varepsilon^2/6\right)^t.$$

Следовательно,

$$E\left[\Phi(t)\right] \le n\exp\left(t\varepsilon^2/6\right) \text{ и } E\left[\Phi(t^*)\right] \le n^{3/2+r/2}.$$

Отсюда по неравенству Маркова

$$\forall \text{ с.в. } \xi \ge 0, t > 0 \to P(\xi \ge t) \le E\xi/t,$$

имеем

$$P\left(\Phi(t^*) \le n^{3/2(1+r)}\right) = 1 - P\left(\Phi(t^*) \ge n^{3/2(1+r)}\right) \ge 1 - E\left[\Phi(t^*)\right]\Big/ n^{3/2(1+r)} \ge 1-\sigma.$$

Тогда, логарифмируя обе части неравенства

$$\exp\left(\varepsilon U_i\left(t^*\right)/2\right) = P_i\left(t^*\right) \le \sum_{i=1}^{n} P_i\left(t^*\right) = \Phi\left(t^*\right) \le n^{3/2(1+r)},$$

имеющего место с вероятностью не меньшей, чем $1-\sigma$, получим

$$P\left(\varepsilon U_i\left(t^*\right)/2 \le 3/2(1+r)\ln n,\ i=1,\dots,n\right) \ge 1-\sigma,$$

$$P\left(U_i\left(t^*\right)/\left(3(1+r)\varepsilon^{-2}\ln n\right) \le \varepsilon,\ i=1,\dots,n\right) \ge 1-\sigma.$$

Откуда уже следует то, что требуется

$$P\left(\mathrm{A}\vec{x}\left(t^*\right) \le \varepsilon\vec{e}\right) \ge 1-\sigma.$$

Авторы выражают благодарность студентам МФТИ Ивану Коноваленко, обнаружившего важную опечатку в работе, и Егору Клочкову, существенно способствовавшего написанию п. 4, а также В.В. Вьюгину, В.В. Высоцкому, А.В. Назину, Ю.Е. Нестерову, А.С. Немировскому, Б.Т. Поляку, А.М. Райгородскому, В.Г. Спокойному, С.П. Тарасову, Р. Темпо, А.А. Трембе, П.Ю. Чеботареву за ряд ценных замечаний. Также хотелось бы поблагодарить студентку ФУПМ МФТИ Айсу Имееву, которая по просьбе авторов провела ряд численных экспериментов с рассматриваемыми в статье алгоритмами.

К статье имеется мини-курс лекций (с видео), который один из авторов прочитал в июле 2013 года в ЛШСМ-2013:

http://www.mathnet.ru/php/presentation.phtml?option_lang=rus&presentid=7259

## Список литературы